\documentclass[11pt]{article}
\usepackage[small]{titlesec}
\usepackage[cm]{fullpage}
\usepackage{amsmath, amssymb, amsfonts}
\usepackage{amscd}
\usepackage{amsmath}
\usepackage{amsfonts}

\newcommand{\de}{\partial}

\newcommand{\ve}{\varepsilon}
\newcommand{\vp}{\varphi}
\newcommand{\mn}{\sqrt{-1}}

\newcommand{\ddbar}{\frac{\sqrt{-1}}{\pi} \partial \overline{\partial}}

\newcommand{\ti}[1]{\tilde{#1}}

\renewcommand{\leq}{\leqslant}
\renewcommand{\geq}{\geqslant}
\renewcommand{\le}{\leqslant}
\renewcommand{\ge}{\geqslant}
\begin{document}
\newcounter{remark}
\newcounter{theor}
\setcounter{remark}{0}
\setcounter{theor}{1}
\newtheorem{claim}{Claim}
\newtheorem{theorem}{Theorem}[section]
\newtheorem{proposition}{Proposition}[section]
\newtheorem{conjecture}{Conjecture}[section]
\newtheorem{lemma}{Lemma}[section]
\newtheorem{defn}{Definition}[theor]
\newtheorem{corollary}{Corollary}[section]
\newenvironment{proof}[1][Proof]{\begin{trivlist}
\item[\hskip \labelsep {\bfseries #1}]}{\end{trivlist}}
\newenvironment{remark}[1][Remark]{\addtocounter{remark}{1} \begin{trivlist}
\item[\hskip
\labelsep {\bfseries #1  \thesection.\theremark}]}{\end{trivlist}}
\setlength{\arraycolsep}{2pt}
\centerline{\bf \large  Plurisubharmonic functions and nef classes }
\smallskip

\centerline{\bf \large on complex manifolds\footnote{Research supported in part by National Science Foundation grants DMS-08-48193 and DMS-10-05457.  The second-named author is also supported in part by a Sloan Foundation fellowship.}}

\bigskip
\bigskip
\centerline{\bf Valentino Tosatti and Ben Weinkove}

\bigskip

\begin{abstract}  We prove the existence of plurisubharmonic functions with prescribed logarithmic singularities on complex 3-folds equipped with a nef class of positive volume.  We prove the same result for rational classes on Moishezon $n$-folds.
\end{abstract}

\section{Introduction}

\bigskip

We recall that a function $f: \Omega \rightarrow [-\infty, \infty)$, for a domain $\Omega \subset \mathbb{C}^n$ is plurisubharmonic if it is upper semi-continuous and for every $a, b \in \mathbb{C}^n$, the map 
$$z \in \mathbb{C} \mapsto f(a+ bz) \in [-\infty, \infty)$$ is subharmonic where it is defined.

Now let $(X, \omega)$ be a compact K\"ahler manifold of complex dimension $n$.   There are complex coordinate charts $B_i$ on which we can write $\omega = \ddbar g_i$ for smooth potential functions $g_i$ (unique up to adding pluriharmonic functions). 
We define a  function $\varphi:X\to [-\infty,+\infty)$  to be $\omega$-plurisubharmonic if $\varphi$ is upper semi-continuous, not identically $-\infty$,  and each $g_i + \varphi$ is plurisubharmonic (clearly, this does not depend on the choice of potentials $g_i$).  A basic fact is that an $L^1$ function $\varphi$ on $X$ satisfying
\begin{equation} \label{pos}
 \omega + \ddbar \varphi \ge 0,
\end{equation}
in the sense of currents, agrees with a unique $\omega$-plurisubharmonic function almost everywhere.  Conversely, every $\omega$-plurisubharmonic $\varphi$  is in $L^1$ and satisfies (\ref{pos}).

The space of $\omega$-plurisubharmonic $\varphi$ has been the focus of considerable study in the last few decades.  Note that every closed positive real $(1,1)$-current cohomologous to $[\omega] \in H^{1,1}(X, \mathbb{R})$ can be written as $\omega+\ddbar \varphi$ for some such $\varphi$.

It was shown in \cite{Ho, Ti} (see also \cite{TY}) that there exists a constant $\alpha>0$ depending only on $[\omega]$ such that
$$\int_X e^{-\alpha \varphi} \omega^n \le C,$$
for all $\omega$-plurisubharmonic $\varphi$.  
This shows in particular that singularities of $\varphi$ can be at most logarithmic.  

An interesting and well-known problem is: 
can we find a $\varphi$  with prescribed logarithmic singularities at given points on $X$?    In the case when $[\omega]$ is the Chern class of a holomorphic line bundle $L$, this is equivalent to prescribing singular Hermitian metrics on $L$.   This can be naturally extended to line bundles $L$ which are only nef.  Deep results in \cite{de, dem2, Si} and others, used the construction of singular $\varphi$ to prove effective results in algebraic geometry.

In this short note we investigate this problem on general complex (non-K\"ahler) manifolds.  Our motivation is to try to understand whether techniques from K\"ahler geometry can be extended to non-K\"ahler complex geometry, at least when natural analogues exist.  

  Suppose now that $X$ is only a compact complex manifold.   Let $\beta$ be a closed real $(1,1)$-form on $X$, consider the (finite-dimensional) real Bott-Chern cohomology group
$$H^{1,1}_{\mathrm{BC}}(X, \mathbb{R})=\frac{\{\beta \textrm{ closed real }(1,1)\textrm{-forms}\}}{\{\beta=\ddbar \psi, \psi\in C^\infty(X,\mathbb{R})\}},$$
and call $[\beta]$ the class of $\beta$ in $H^{1,1}_{\mathrm{BC}}(X, \mathbb{R})$. 

Since a positive  $\beta$ only exists if  $X$ is K\"ahler we consider instead the case when the class 
 $[\beta]$ is nef (as defined in \cite{dp}), which means that 
for any $\ve>0$ there exists
a representative $\beta+\ddbar \psi_\ve$ so that $\beta+\ddbar\psi_\ve > -\ve\omega$, where $\omega$ is some fixed Hermitian metric on $X$.  

The closed form $\beta$ admits local potential functions and thus we can define the notion of $\beta$-plurisubharmonic in the same way as described above.  Under the cohomological assumption that $\int_X \beta^n>0$, 
we look for $\beta$-plurisubharmonic functions
$\varphi$ with prescribed logarithmic singularities.

Our main result is the following:

\bigskip

\pagebreak[3]
\noindent
{\bf Main Theorem} \, \emph{Let $X$ be a compact complex manifold of dimension $n$.   Suppose there exists a class $[\beta] \in H^{1,1}_{\mathrm{BC}}(X, \mathbb{R})$ which is nef and satisfies $\int_X \beta^n>0$.   Assume {\bf either}
\\ \\ \indent
(i)  $n=2$ or $n=3$. \\
 \\
{\bf or}   (ii)  $X$ is Moishezon and $[\beta]  \in H^{1,1}_{\mathrm{BC}}(X, \mathbb{Q}):=
H^{1,1}_{\mathrm{BC}}(X, \mathbb{R})\cap H^2(X,\mathbb{Q})$. \\ \\
Fix points $x_1, \ldots, x_N \in X$ and choose positive real numbers $\tau_1, \ldots, \tau_N$ so that
\begin{equation}\label{condition}
\sum_{j} \tau_j^n < \int_X \beta^n.
\end{equation}
 Then there exists a $\beta$-plurisubharmonic $\varphi$  with logarithmic poles at $x_1, \ldots, x_N$:
\begin{equation}\label{want}
\varphi(z) \le \tau_j \log  |z| + O(1),
\end{equation}
in a coordinate neighborhood $(z_1, \ldots, z_n)$ centered at $x_j$, where $|z|^2 = |z_1|^2+ \cdots + |z_n|^2$. In particular, the Lelong number of $\varphi$ at each point $x_j$ is at least $\tau_j$.}

\bigskip

Recall that a Moishezon manifold is a compact complex manifold which is bimeromorphic to a projective manifold.  An equivalent definition is that a Moishezon manifold is a compact complex manifold admitting a big line bundle $L$ (meaning $\textrm{dim}\,  H^0(X, L^k) > c k^n$ for $k$ large, for some fixed $c>0$).

If $X$ is K\"ahler (without imposing (i) or (ii)), this result is due to Demailly \cite{de} whose proof made use of Yau's solution of the complex Monge-Amp\`ere equation on K\"ahler manifolds \cite{Ya}.  

In the case of one point ($N=1$), inequality \eqref{condition} is sharp in general: taking $X=\mathbb{CP}^n$ and $[\beta]$ to be the (ample) anticanonical class (so $\beta$ is $n+1$ times the Fubini-Study metric), equation \eqref{condition} says that $\tau<n+1$, and it is well known (see e.g. \cite[Proposition 2.1]{CG}) that $n+1$ is indeed the maximum order of logarithmic pole of any $\beta$-plurisubharmonic function.

We remark that our result in case (i) is only really new in the case $n=3$.   The reason is that in dimension $2$, any surface $X$ as in the Main Theorem is necessarily K\"ahler. In fact, the existence of a closed real $(1,1)$-form $\beta$ with $\int_X\beta^2>0$ implies that the intersection form on $H^{1,1}(X, \mathbb{R})$ is not negative definite, and thanks to a classical theorem of Kodaira \cite{Ko} this implies that $b_1(X)$ is even. A theorem of Miyaoka-Siu \cite{M, Si2} (see also \cite{Bu, La}) then implies that $X$ is K\"ahler.

On the other hand there are certainly many non-K\"ahler $3$-folds satisfying the hypotheses of the Main Theorem  (see for example \cite{Hi} and the description in \cite[Example 3.4.1, p.443]{Ha}).  Indeed, such manifolds were discussed by Demailly-P\u{a}un  in \cite{dp}, where it was   conjectured  that a compact complex $n$-fold $X$ with a nef class $[\beta]$ of type $(1,1)$ with $\int_X \beta^n >0$ is bimeromorphic to a K\"ahler manifold.

Demailly \cite{de} proved the K\"ahler version of the Main Theorem using Yau's existence result for solutions to the complex Monge-Amp\`ere equation on K\"ahler manifolds  \cite{Ya}.  In this paper we follow along the same lines of argument as Demailly, but now apply the recent extension of Yau's Theorem to general complex manifolds \cite{tw} (see also \cite{Ch, GL, tw2, dk, bl, Gi}).  However, a difficulty arises here in the non-K\"ahler case due to the  fact that for a $(1,1)$-form $\Omega$ and a function $f$, the equality
$$\int_X \left(\Omega + \ddbar f\right)^n = \int_X \Omega^n$$
does not hold in general if $\Omega$ is not closed.  We can overcome this obstacle in dimensions 2 and 3 by making use of Gauduchon metrics.

We conjecture that the Main Theorem holds for any dimension without the Moishezon and rationality assumptions in (ii).

\section{Proof of the Main Theorem}

Let $\omega$ be a Gauduchon metric on $X$, which means that $\partial \overline{\partial} (\omega^{n-1})=0$ (such a metric always exists \cite{Ga}).  Following Demailly \cite{de}, we  choose coordinates $z^1, \ldots, z^n$ in a neighborhood centered at $x_j$ and put
$$\gamma_{j,\ve} = \frac{\mn}{\pi}\de\overline{\de} \left(\chi\left(\log\frac{|z|}{\ve}\right) \right),$$
where $|z|^2 = |z^1|^2 + \cdots + |z^n|^2$, and $\chi:\mathbb{R} \rightarrow \mathbb{R}$ is smooth, convex, increasing and satisfies $\chi(t)=t$ for $t \ge 0$ and $\chi(t)=-1/2$ for $t \le -1$.  Then observe that $\gamma_{j,\ve}=0$ if $|z| > \ve$, so we can extend it to zero on the whole of $X$. This way, $\gamma_{j,\ve}$ is a closed nonnegative smooth $(1,1)$-form on $X$ that satisfies
$$\int_X\gamma_{j,\ve}^n= \int_{|z| \le \ve} \gamma_{j,\ve}^n = 1,$$
and $\gamma_{j,\ve}^n\rightharpoonup \delta_{x_j}$ as $\ve\to 0$.
Since $\beta$ is nef, for any $\ve>0$ there exists a smooth function $\psi_\ve$ such that $\beta + \ve \omega+\ddbar\psi_\ve$ is Hermitian.  

Now recall that the Hermitian version of Yau's theorem \cite{tw} states that given a Hermitian metric $\hat{\omega}$ and a smooth function $F$ there exists a unique smooth $f$ and a unique constant $K>0$ solving
$$\left(\hat{\omega}+ \ddbar f\right)^n = K e^F \hat{\omega}^n, \quad  \hat{\omega}+ \ddbar f>0, \quad \sup_X f=0.$$ 
Note that, in the case of Yau's theorem where $\hat{\omega}$ is K\"ahler, integration by parts shows that  $K = \int_X \hat{\omega}^n/ \int_X e^F \hat{\omega}^n$.  In the non-K\"ahler case, no such formula holds in general and this is the source of the difficulty.
 
Applying this with reference metric $\beta + \ve \omega+\ddbar\psi_\ve$, we obtain a smooth  $\varphi_{\ve}$ with
$$\left( \beta + \ve \omega + \ddbar (\psi_\ve+\varphi_{\ve})\right)^n = C_{\ve} \left( \sum_j {\tau_j^n} \gamma_{j,\ve}^n + \delta \omega^n\right)$$
and
$$\beta + \ve \omega + \ddbar (\psi_\ve+\varphi_{\ve}) >0,  \quad \sup_X (\psi_\ve+\varphi_{\ve})=0,$$
where $\delta>0$ is fixed and $C_{\ve}$ is a (uniquely determined) positive constant.   The key fact that we need now is the lower bound $C_{\ve} \ge 1$ for $\ve$ and $\delta$ sufficiently small.  We remark that there are some formal similarities between the argument given here for this lower bound and the proofs of Proposition 3.8 in \cite{dps} and Theorem 1.2 in \cite{J}.

We consider first the case (i) when $n=2$. Then $\de\overline{\de}\omega=0$ and, using  the fact that $\beta$ is closed, we have
$$C_\ve=\frac{\int_X( \beta + \ve \omega + \ddbar (\psi_\ve+\varphi_{\ve}))^2}{\int_X  \sum_j {\tau_j^2} \gamma_{j,\ve}^2 + \delta \omega^2}=\frac{\int_X( \beta + \ve \omega)^2}{\sum_j {\tau_j^2} + \delta \int_X \omega^2}\geq\frac{\int_X\beta^2 - \ve^2 \int_X \omega^2}{\sum_j {\tau_j^2} + \delta \int_X \omega^2},$$
since
$$\int_X (\beta + \ve \omega)^2 = \int_X \beta^2 + 2 \ve \int_X \left(\beta + \ve \omega + \ddbar \psi_{\ve} \right) \wedge \omega - \ve^2 \int_X \omega^2.$$
Choosing $\delta$ sufficiently small and using \eqref{condition}, we get for any $\ve>0$ sufficiently small,
$$C_{\ve} \geq 1.$$

We now consider the case when $n=3$. For simplicity call $\beta_\ve= \beta + \ve \omega + \ddbar (\psi_\ve+\varphi_{\ve})$, and note that
$$\beta+\ddbar (\psi_\ve+\varphi_{\ve})=\beta_\ve -\ve\omega,$$
and so 
$$\int_X(\beta_\ve -\ve\omega)^3=\int_X \left(\beta+\ddbar (\psi_\ve+\varphi_{\ve})\right)^3=\int_X\beta^3>0.$$
Using the Gauduchon condition $\de\overline{\de}(\omega^2)=0$, we see that
\begin{equation}\label{est1}
\int_X\beta_\ve\wedge\omega^2=\int_X (\beta+\ve\omega)\wedge\omega^2
\leq C\int_X\omega^3.
\end{equation}
On the other hand, using \eqref{est1}, we have
\begin{eqnarray} \nonumber
\int_X\beta_\ve^3&=&\int_X(\beta_\ve-\ve\omega+\ve\omega)^3\\ \nonumber
&=& \int_X(\beta_\ve-\ve\omega)^3+3\ve\int_X\beta_\ve^2\wedge\omega
-3\ve^2\int_X\beta_\ve\wedge\omega^2+\ve^3\int_X\omega^3\\ \nonumber
&\geq&\int_X \beta^3-C'\ve^2,
\end{eqnarray}
and so
$$C_{\ve}=\frac{\int_X\beta_\ve^3}{\int_X  \sum_j {\tau_j^3} \gamma_{j,\ve}^3 + \delta \omega^3} \ge   \frac{\int_X \beta^3 - C' \ve^2}{ \sum_j {\tau_j^3} + \delta \int_X \omega^3},$$
and hence  choosing $\delta$ sufficiently small and using \eqref{condition}, we get for $\ve>0$ sufficiently small,
$$C_{\ve} \ge 1.$$

Finally we prove $C_{\ve} \ge 1$ in the case (ii).   Because $[\beta]$ is rational we have $\ell [\beta] = c_1(L)$ for some line bundle $L$ over $X$ and some integer $\ell\geq 1$.  In this case there exists a Hermitian metric $h_{\ve}$ on $L$ such that $\frac{1}{\ell} c_1(L, h_{\ve}) = \beta + \ddbar (\psi_{\ve} + \varphi_{\ve})$, where $c_1(L, h_{\ve})$ is the curvature form of the Hermitian metric $h_{\ve}$.  Denote by $X(0)$ the set of $x \in X$ such that $c_1(L, h_{\ve})$ has $0$ negative eigenvalues.   We now apply Demailly's holomorphic Morse inequalities \cite{de0} to see that for $k$ large we have
\begin{eqnarray}\label{mor} \nonumber
\dim H^0(X,L^k) & \le & \frac{k^n}{n!}\int_{X(0)} c_1(L, h_{\ve})^n+o(k^n) \\ \nonumber
& \le & \frac{k^n}{n!} \int_{X(0)} \left(c_1(L, h_{\ve}) + \ve \ell\omega \right)^n +o(k^n)\\ \nonumber
&  =  & \frac{C_{\ve}\ell^n k^n}{n!} \int_{X(0)} \left( \sum_j {\tau_j^n} \gamma_{j,\ve}^n + \delta \omega^n\right)+o(k^n) \\
& \le & \frac{C_{\ve}\ell^n k^n}{n!} \left( \sum_j \tau_j^n +  \delta \int_X \omega^n \right)+o(k^n).
\end{eqnarray}
We now estimate the number of sections of $L^k$ using the Riemann-Roch theorem. 
Since the manifold $X$ is Moishezon there exists a modification $\mu:\ti{X}\to X$ with $\ti{X}$ a projective manifold. The pullback $\mu^*L$ is then a nef line bundle on $\ti{X}$ with $\int_{\ti{X}}c_1(\mu^*L)^n=\ell^n\int_X\beta^n>0$.
Because $\mu^*L$ is nef, its higher cohomology groups satisfy
$\dim H^q(\ti{X},\mu^*L^k)=O(k^{n-1}),$ for $q>0$ (see Example 1.2.36 in \cite{Laz}), and a standard Leray spectral sequence argument (see (2.1) in \cite{ag}) 
shows that $\dim H^q(X,L^k)=O(k^{n-1}),$ for $q>0$.
By the Riemann-Roch theorem we now have
\begin{equation}\label{rr}
\dim H^0(X,L^k)=\frac{\ell^n k^n}{n!}\int_X\beta^n+o(k^n).
\end{equation}
Combining \eqref{mor} and \eqref{rr} and taking $k$ large we get
$$\int_X\beta^n\leq C_{\ve}\left( \sum_j \tau_j^n +  \delta \int_X \omega^n \right).$$
Choosing $\delta>0$ sufficiently small we obtain $C_{\ve} \ge 1$.

We claim that $\psi_\ve+\varphi_{\ve}$ is uniformly bounded in $L^1$. This is because we can choose a large constant $A$ so that $A\omega\geq \beta+\ve\omega$
for all $0<\ve\leq 1$, and then the function $\psi_\ve+\varphi_{\ve}$ satisfies
$A\omega+\ddbar (\psi_\ve+\varphi_{\ve})>0$, and then Proposition 2.1 in \cite{dk} (for example) gives a uniform $L^1$ bound that depends only on $X, A, \omega$.

 This implies that there is a sequence $\ve_k \rightarrow 0$ such that $\psi_{\ve_k}+\varphi_{\ve_k}$ converges in $L^1$ to a $\beta$-plurisubharmonic function $\varphi$.  Indeed, we can recover this from the local statement about compactness of plurisubharmonic functions in a domain in $\mathbb{C}^n$ which are uniformly bounded in $L^1$, in the following way: we cover $X$ with finitely many coordinate charts $B_i$ so that on each $B_i$ there is a smooth function $\rho_i$ with
$\ddbar\rho_i > \beta+\omega$. Then on each $B_i$ the functions $\rho_i+\psi_{\ve}+\varphi_{\ve}$ are plurisubharmonic and uniformly bounded in $L^1$ (independent of $\ve$), so the local statement applies. 
  
We now show that $\varphi$ has the desired logarithmic singularities.

Take $\Omega$ a neighborhood of $x_j$ (which we can assume contains the set $\{|z|<1\}$) and consider the smooth plurisubharmonic function on $\Omega$
$$u = C_{\ve}^{1/n} \tau_j (\chi (\log(|z|/\ve)) + \log \ve) + C_1$$
for $C_1$ a large constant.  Let $h$ be a smooth function on $\overline{\Omega}$ with
$$\ddbar h \ge  \beta+ \omega$$
and put $v= h+ \psi_\ve+\varphi_{\ve}$, so $v$ is a smooth plurisubharmonic function on $\Omega$ which is bounded from above by $C_0$, say.
 
 Then if $\ve$ is sufficiently small and $C_1$ is sufficiently large we have 
$$u|_{\partial \Omega}=C_\ve^{1/n}\tau_j \log|z|+C_1 \geq C_0\geq v|_{\partial \Omega}.$$
But in addition we have on $\Omega$,
$$\left(\ddbar v\right)^n \ge \left(\beta + \ve \omega + \ddbar (\psi_\ve+\varphi_{\ve})\right)^n \ge C_{\ve} \tau_j^n \gamma_{j, \ve}^n = \left(\ddbar u\right)^n.$$
Then by the Bedford-Taylor comparison principle for Monge-Amp\`ere (e.g. Lemma 6.7 in \cite{de}),
$$u \ge v \quad \textrm{on} \ \Omega,$$
and hence when $|z|<1/2$ and $\ve$ is small we have
$$\psi_\ve+\varphi_{\ve} \leq C_{\ve}^{1/n} \tau_j (\chi (\log(|z|/\ve)) + \log \ve)
+C_2 \le C_{\ve}^{1/n} \tau_j \log (|z|+\ve) + C_2 \le \tau_j  \log (|z|+\ve) +C_2,$$
where we are using the fact that $C_{\ve} \ge 1$. Since $\psi_\ve+\varphi_{\ve}$ converges to $\vp$ in $L^1$, Hartogs' Lemma (see for example Proposition 2.6 (2) in \cite{GZ})
implies \eqref{want}.
Q.E.D.\\

Let us remark here that, in any dimension, the constant $C_\ve$ is always bounded above (independent of $\ve$ but depending on $\delta$) when $\ve$ is small. Indeed, at the point on $X$ where the function $\psi_\ve+\varphi_{\ve}$ achieves its maximum, we have that $\beta+\ve\omega$ is positive definite and moreover
$\ddbar(\psi_\ve+\varphi_{\ve})\leq 0$ and so at that point we have
$$C\omega^n \geq (\beta+\ve\omega)^n\geq \left(\beta+\ve\omega+\ddbar(\psi_\ve+\varphi_{\ve})\right)^n=
C_\ve \left( \sum_j {\tau_j^n} \gamma_{j,\ve}^n + \delta \omega^n\right)\geq
C_\ve \delta\omega^n,$$
giving $C_\ve\leq C/\delta$.

A similar argument only gives a lower bound of the form $C_\ve\geq c\ve^{3n}$, in any dimension. In fact, a direct calculation shows that on $X$ we have
$$\gamma_{j,\ve}^n\leq \frac{C}{\ve^{2n}}\omega^n,$$
for a constant $C$ independent of $\ve$. We then compute
$$\left(\beta+\ve\omega+\ddbar(\psi_\ve+\varphi_{\ve})\right)^n=C_\ve\left( \sum_j {\tau_j^n} \gamma_{j,\ve}^n + \delta \omega^n\right)\leq
\frac{C}{\ve^{2n}} C_\ve \omega^n.$$
Recall that since $[\beta]$ is nef there is a smooth function $\psi_{\ve/2}$
so that $\beta+\frac{\ve}{2}\omega+\ddbar\psi_{\ve/2}>0$. We then write
$$\beta+\ve\omega+\ddbar(\psi_\ve+\varphi_{\ve})=
\left(\beta+\frac{\ve}{2}\omega+\ddbar\psi_{\ve/2}\right)+\left(\frac{\ve}{2}\omega+
\ddbar(\psi_\ve+\varphi_{\ve}-\psi_{\ve/2})\right).$$
At the point on $X$ where $\psi_\ve+\varphi_{\ve}-\psi_{\ve/2}$ achieves its minimum we have $\ddbar(\psi_\ve+\varphi_{\ve}-\psi_{\ve/2})\geq 0$, and so at that point
$$\frac{\ve^n}{2^n}\omega^n\leq\left(\beta+\ve\omega+\ddbar(\psi_\ve+\varphi_{\ve})\right)^n
\leq \frac{C}{\ve^{2n}} C_\ve \omega^n,$$
which gives $C_\ve\geq c\ve^{3n}$.

\bigskip
\noindent
Mathematics Department, Columbia University, 2990 Broadway, New York, NY 10027

\bigskip
\noindent
Mathematics Department, University of California at San Diego, 9500 Gilman Drive, La Jolla, CA 92093

\end{document}